\documentclass[12pt]{article}

\usepackage[a4paper]{geometry}
\usepackage{ajc-FR}
\usepackage{amsmath,amsfonts,latexsym}
\usepackage{graphicx}
\usepackage{placeins}

\Volume{XY(Z)}
\Year{2019}
\firstpage{1}
\lastpage{8}
\Revised{3 December 2019}
\runninghead{\copyright \, F. Rowley 2019\,/\,FRRGC DRAFT arXiv version, pages 1 -- 8}
\numberwithin{equation}{section}

\newtheorem{theorem}{Theorem}[section]

\newenvironment*{proof}
{\begin{list}{}{\setlength{\leftmargin}{0em}\setlength{\rightmargin}{0em}}
\item[] {\sc Proof:}} {\hfill$\Box$
\end{list}}

\begin{document}

\title{\bf Some further results in \\Ramsey graph construction}
\author{Fred Rowley\thanks{formerly of Lincoln College, Oxford, UK.}}
\Addr{West Pennant Hills, \\NSW, Australia.\\{\tt fred.rowley@ozemail.com.au}}


\date{\dateline{18 Feb 2019}{DD Mmm CCYY}}

\maketitle


\begin{abstract}
A construction described by the current author (2017) uses two linear prototypes to build a compound graph with Ramsey properties inherited from the prototype graphs.  The resulting graph is linear; and cyclic if both prototypes are cyclic.  However, it will not generate a cyclic graph from a general linear prototype.  

Building on the properties of that construction, this paper proves that a general linear prototype graph of order $m$ can be extended using a single new colour to produce a new cyclic graph of order $3m - 1$ which is triangle-free in the new colour, and has the same clique-number as the prototype in every other colour. 

The paper then describes a cyclic Ramsey $(3,3,4,4; 173)$-graph derived by constrained tree search -- thus proving that $R(3,3,4,4) \ge 174$.  Using a quadrupling construction to produce a further cyclic graph, it is shown that $R(3,4,5,5) \ge 693$.  

A compound cyclic Ramsey $(3,7,7; 622)$-graph derived by a limited manual search is then described.  Further construction steps produce a $(8,8,8; 6131)$-graph, showing that $R_3(8) \ge 6132$.

The paper concludes by showing that $R_4(7) \ge 81206$ and $R_4(9) \ge 630566$, implying corresponding improvements in the lower bounds for $R_5(7)$ and $R_5(9)$ and beyond.  These results follow from the existence of cyclic prototype graphs derived by Mathon-Shearer 'doubling'.

\end{abstract}
\bigskip
\small DRAFT \copyright Fred Rowley  February 2019.
\bigskip

\section{Introduction}
This paper addresses the properties of undirected loopless graphs with edge-colourings in an arbitrary number of colours, and the corresponding multicolour classical Ramsey numbers.  

The construction described in \cite{Rowley1} allows the creation of Ramsey graphs with specific properties by combining the distance sets of two linear prototypes.  The resulting graph is linear; and cyclic if both prototypes are cyclic.  However, it will not generate a cyclic graph from a general linear prototype -- a property that is useful on occasions if subsequent constructions require cyclic prototypes.  This paper describes a development of the previous construction which does so in a significant special case.  

It also records some further computational findings and results in the construction of linear and cyclic graphs. 

Notation is defined in section 2.  

In section 3, starting from the construction in \cite{Rowley1}, it is proved that a general linear prototype graph can be extended using a single new colour to produce a new {\bf cyclic} graph of order $3m - 1$ which is triangle-free in the new colour, and has the same clique-number as the prototype in every other colour. This result is of some practical use in providing graphs with known characteristics for use as prototypes in further constructions.  Its distance profile, and variations on it, provide a broad theme for the paper.  

In section 4, the distance sets of a cyclic Ramsey $(3,3,4,4; 173)$-graph are listed -- thus proving that $R(3,3,4,4) \ge 174$.  These sets were found by a non-exhaustive tree search, subject to a heuristically-derived constraint, which made search times manageable.  The form of the constraint was inspired by the broad characteristics of colourings related to the construction described in section 3.  It is further deduced that $R(3,4,5,5) \ge 693$ using a cyclic 'quadrupling' technique related to that described in \cite{XXER}.  

In section 5, the existence a Ramsey $(8,8,8; 6131)$-graph is established.  We start from a prototype $(3,7,7)$-graph, use the quadrupling construction twice, and then delete vertices.  Thus it is proved that $R_3(8) \ge 6132$.

In section 6, the key Tables from \cite{Rowley1} are updated to reflect further work on obtaining linear prototype graphs.  It is proved that $R_4(7) \ge 81206$ and $R_4(9) \ge 630566$, noting corresponding improvements in the lower bounds for $R_5(7)$ and $R_5(9)$ and beyond.  These results follow from the existence (in all cases) of cyclic prototype subgraphs of the graphs obtained by the most basic Mathon-Shearer 'doubling' construction.  

\section{Notation}
In this paper, 

$K_n$ denotes the complete graph with order $n$.  

If $U$ denotes a complete graph with $m$ vertices $\{u_0, {\dots} , u_{m-1}\}$, then: 

A {\it (q-)colouring} of $U$ is a mapping of the edges $({u_i}, {u_j})$ of $U$ into the set of integers $s$ where $1\le s \le q$.  

The {\it distance} between two vertices ${u_i}, {u_j}$, or, equivalently, the {\it length} of the edge $({u_i}, {u_j})$ connecting them, is defined as $\mid j - i \mid$.

A colouring of $U$ is {\it linear} if and only if the colour of any edge $({u_i}, {u_j})$ depends only on the length of that edge.  In such a case the colour of an edge of length $l$ is written $c(l)$.  

A colouring of $U$ is {\it cyclic} if and only if (a) it is linear, and (b) $c(l)=c(m-l)$ for all l such that $1 \le l \le m-1$.  

The {\it clique number} of graph $U$ in colour $s$ is the largest integer $i$ such that $U$ contains a subgraph which is a copy of $K_i$ in that colour.  

A {\it Ramsey graph} $U(k_1, {\dots} \, ,k_r; m)$, with all $k_s \ge 2$, is a complete graph of order $m$ with a colouring such that for each colour $s$, where $1 \leq s \leq r$, there exists no complete monochrome subgraph $K_q$ of $U$ in the colour $s$ for any $q \geq k_s$.  Equivalently, the clique number of $U$, in any colour $s$, is strictly less than $k_s$.  Such a graph $U$ may conveniently be described as a $(k_1, {\dots} , k_r; m)${\it-graph}.  

The {\it Ramsey number} $R(k_1, \dots \, , k_r)$ is the unique lowest integer $m$ such that no \newline
 $U(k_1, {\dots} \, , k_r; m)$ exists.

\section{Construction of Cyclic Graphs from Linear Graphs}

\begin{theorem}
  \label{Thm:C-Thm2}
(Construction Theorem)\\
Given any linear Ramsey graph $U(k_1, k_2, {\dots} , k_r; m)$, it is possible to construct a cyclic Ramsey graph $W(k_1, k_2, {\dots} , k_r, 3; 3m-1)$.
\end{theorem}

The theorem depends on a relatively simple construction process which adds $2m-1$ vertices and their incident edges, including a single new colour, to the linear prototype graph $U$.   

\begin{proof}
We begin by considering the set of lengths of all the edges of $U$, which we call $L$, consisting of the integers $\{\, l \mid 1 \le l \le m-1 \,\}$.  A linear colouring gives rise to a natural partition of that set into subsets $L_s$ containing the lengths of edges of each colour $s$. That is, for $1 \leq s \leq r \,$:

\hspace{120pt} $L_s = \{\, l \mid c(l) = s \}.$

It is a well-known result that any linear graph $U$ contains a copy of $K_s$ in colour $s$ if and only if there exists a subset of the set $L_s$ of order $s-1$ such that each of the members of the subset and all of their non-zero pairwise differences are contained in $L_s$.  For if such a subset exists, one can construct a set of all the vertices $u_i \in U$ having index-numbers $i$ in the subset.  Taking the union of that set of vertices with $u_0$ gives us the vertices of a copy of $K_{k_r}$ in $U$.  The converse is essentially proved by reversing the process, having first selected (using linearity) a copy of $K_{k_s}$ with a vertex set that includes $u_0$.  

This result provides the basis for our proof.

Using the construction in \cite{Rowley1} we first construct a linear $(k_1, k_2, {\dots} , k_r, 3; 3m-1)$-graph $V$ with vertices $v_i$, for $0 \le i \le 3m-2$.  The set of lengths of all the edges of $V$ may be called $L'$ and consists of the integers $\{\, l \mid 1 \le l \le (3m-2) \,\}$, which are partitioned into distinct subsets according to colour, as follows:

We define a subset of $L'$ which we call $L'_{r+1} = \{ l \mid m \le l \le (2m-1) \}$.   

We further define subsets of $L'$, for each $s$, where $1 \le s \le r$, as follows:

\hspace{120pt} $A'_s = \{\, l \mid l \in L_s\}$

\hspace{120pt} $B'_s = \{\, l + (2m-1) \mid l \in L_s \}$, and

\hspace{120pt} $L'_s = A'_s \cup B'_s$.	

It is easy to verify that $\bigcup L'_t = L'$, where the union includes all colours.

The new graph $V$ is clearly linear. From the proof in \cite{Rowley1} it follows that the clique number of $V$ in any colour $s$ where $1 \le s \le r$ is the same as for $U$.  For colour $r+1$ the clique number is clearly $2$, since no two members of $L'_{r+1}$ can have a difference in $L'_{r+1}$.  

Now we define a second graph $W$, also with $3m-1$ vertices $w_i$.  The set of lengths of all the edges of $W$ may be called $L''$.

We define the subset $L''_{r+1} = L'_{r+1}$.   

We further define, for each $s$, where $1 \le s \le r$, 

\hspace{120pt} $A''_s = A'_s $,

\hspace{120pt} $B''_s = \{\, (5m-2) - l \mid l \in B'_s \}$, and

\hspace{120pt} $L''_s = A''_s \cup B''_s$.

The new graph $W$ is clearly well-defined and linear. Because $B''_s$ contains the complements of members of $B'_s$ with respect to $5m-2$, $W$ is also cyclic.  

We have assumed that there is no subset of $U$ that is a monochromatic copy of $K_{k_s}$ in colour $s$. We aim to prove there can be no such subset in $W$.  

Assume to the contrary that there is such a copy ($H$, say) in $W$.  Consider the set of index-numbers of its vertex-set, $\{j_1, j_2, \dots, j_{k_s}\}$.  

If any pair of these index-numbers has an absolute difference greater than $m-1$, then that difference must be at least $2m$.  If no such pair exists, then the lengths of all the edges of $H$ must be less than $m$ and therefore there must be an identical copy of $K_s$ in the same colour within $U$ (with the same set of index-numbers): which is a contradiction.  

Therefore $H$ must have at least one edge of length at least $2m$.  If so, then we can partition the set of index-numbers of $H$ into two non-empty subsets $S_1 = \{j_1, \dots, j_p \}$ and $S_2 = \{j_{p+1}, \dots, j_{k_s} \}$.  This partition is made on the basis that the length of an edge joining any member of $S_1$ to any member of $S_2$ is of length at least $2m$.  It is a straightforward consequence of the colouring that there can be no more than two such subsets.  We may assume without loss of generality that the index-numbers of these subsets are strictly increasing.  

We now define a mapping from the vertices of $W$ to the vertices of $V$ as follows:

For $1 \le t \le p$, define $w_{j_t} \rightarrow v_{j_p - j_t}$. Thus if there is an edge in $W$ joining $w_{j_x}$ and $w_{j_y}$, both index-numbers being members of $S_1$, then the length of the edge in $V$ that joins their images is the same as in $W$, and the colours of those two edges are both $s$.  

For $p+1 \le t \le s$, define $w_{j_t} \rightarrow v_{((5m-2)-(j_t - j_p))}$. Again, if there is an edge in $W$ joining $w_{j_x}$ and $w_{j_y}$, both index-numbers being members of $S_2$, then the length of the edge in $V$ that joins their images is the same as in $W$, and the colours of those two edges are both $s$.  

The key remaining issue is how this mapping transforms the lengths of images of edges joining, say, $w_{j_x}$ and $w_{j_y}$, where $j_x \in S_1$ and $j_y \in S_2$.

In that case, the length of the image in $V$ of this edge is $(5m-2)-(j_y - j_p)-(j_p - j_x) = (5m-2) - (j_y-j_x)$.  We can see that $2m \le (j_y-j_x) \le 3m-2$. From the method of construction of $B''$ we know that this means that the colours of the edge and its image in $V$ are again both $s$.

Thus the image of $H$ in $V$ is a copy of $K_{k_s}$ in colour $s$, which is another contradiction.  This completes the proof.  

\end{proof}

This quite simple result has been of some practical use in providing cyclic graphs with known characteristics for use as prototypes in further constructions.  The overall profile of the distances -- with one colour concentrated in a mid-range, sometimes also featuring some outlying distances -- is a theme of all the constructions featured in this paper.  The existence of outliers provides scope for defining broader search spaces, as illustrated in sections $4$ and $5$.  

\section{New Lower Bounds for $R(3,3,4,4)$ and $R(3,4,5,5)$}
	
The distance sets for a cyclic $(3,3,4,4; 173)$-graph derived by the author are listed below.  Colours for distances greater than $86$ are implied by the symmetry. The implied lower bound of 174 for $R(3,3,4,4)$ exceeds the current best lower bound quoted in \cite{RadzDS}.  

\begin{table}[!ht]
Colour 1:

{\tt ~2, ~6, ~9, 10, 17, 21, 24, 25, 28, 32,}\\{\tt 39, 40, 55, 62, 75} \medskip
                                                                                
Colour 2:

{\tt 49, 56, 59, 63, 64, 66, 67, 69, 70, 71,}\\{\tt 72, 73, 74, 76, 77, 78, 79, 80, 81, 82,}\\{\tt 83, 84, 85, 86} \medskip
                                                                                
Colour 3:

{\tt ~1, ~5, 11, 12, 15, 19, 20, 22, 27, 29,}\\{\tt 30, 34, 37, 38, 44, 48, 50, 51, 54, 58,}\\{\tt 60, 61, 68} \bigskip
                                                                                
Colour 4:

{\tt ~3, ~4, ~7, ~8, 13, 14, 16, 18, 23, 26,}\\{\tt 31, 33, 35, 36, 41, 42, 43, 45, 46, 47,}\\{\tt 52, 53, 57, 65}       

  \caption{\label{Table:01} Distance Sets for a Cyclic (3,3,4,4;173)-graph.}
\end{table}

\FloatBarrier

This graph was one of several discovered through a non-exhaustive constrained tree search.  The constraint imposed was merely that colour $2$ (one of the triangle-avoiding colours) cannot be used for a distance less than $48$.  

This simple heuristic constraint was inspired by the observation that in generating triangle-free graphs, it is often a feature (loosely stated) that a cluster of distances in the mid-range are of a single common colour, and that the lesser and greater distances largely avoid that colour.  Previous papers featuring distance-based searches (notably \cite{FrSw}) have demonstrated this feature.  

Searches with lower exclusion thresholds have not yet yielded better lower bounds.  Exhaustive searches have not been possible because of the resulting time constraints.

By application of a cyclic 'quadrupling' construction closely related to that featured in Corollary 3 in \cite{XXER}, a $(3,4,5,5; 692)$-graph can be constructed, proving that $R(3,4,5,5) \ge 693$.  
  
\section{A New Lower Bound for $R_3(8)$}
The construction described in section 3 can be usefully varied by modifying the definition of $L''_{r+1}$ in particular cases so that the order of the resulting graph $W$ is increased.

In the case where the prototype is the well-known (7,7; 202)-graph obtained by the basic Mathon-Shearer 'doubling' construction, an interesting result is obtained by defining:

\hspace{40pt} $L''_{r+1} = L''_3 = \{\, 202, 203, \dots , 207, 219, 220, \dots , 403, 415, 416, \dots , 420 \}$.

It can be seen that this retains the general shape of a cluster of distances in the mid-range, as mentioned in section 4, while leaving some gaps in order to include outliers.  By doing this it has so far proved possible to extend the order of the resulting cyclic $(3,7,7)$-graph to $622$ using only a short manual search.  The colours of distances from $1$ to $201$ are the same as those of the $(7,7; 202)$-graph, as are the colours of distances from $421$ to $621$.  The colours of distances from $208$ to $218$ are copies of those from $6$ to $16$ and the colours of distances $404$ to $414$ are a reflection of them.  

Next, we broadly follow the process of Theorem 7 in \cite{XSR2}.  By applying the quadrupling construction featured in Corollary 5 in \cite{XXER} twice in succession, we obtain firstly a $(5,7,8; 2488)$-graph and then a $(9,8,8; 9952)$-graph.  By inspection of the construction at each stage, we can see that the degree of the first vertex in these graphs in colour $3$ is $1214$ and $6131$ respectively. (In fact, the degrees of all vertices are equal in any colour, but that is not necessary to the proof.) These degrees have been validated by computer testing.  Consider the subgraph induced in the latter case on the vertices forming the neighbourhood of the first vertex in colour $3$.  Since $k_1 = 9$ this must be an $(8,8,8; 6131)$-graph, which demonstrates that $R_3(8) \ge 6132$.  This again exceeds the current best lower bound quoted in \cite{RadzDS}. 

\section{Further Results from an Earlier Construction}
The Tables below update the results of the previous paper \cite{Rowley1}, allowing for the inclusion of linear isomorphic images of the well-known $(7,7; 202)$- and $(9,9; 562)$-graphs obtained by the most basic form of the Mathon-Shearer 'doubling' construction. The reader should refer to \cite{Mathon} for the basic construction: the subgraphs can be obtained from all such graphs by a straightforward rearrangement of vertices. 
Numbers revised since the publication of \cite{Rowley1} are shown in blue. As before, bold text indicates numbers exceeding those shown in the Radziszowski Dynamic Survey \cite{RadzDS}.

The inclusion of these graphs allows the proof that $R_4(7) \ge 81206$ and $R_4(9) \ge 630566$ by straightforward application of the methods of \cite{Rowley1}, with correspondingly improved lower bounds for $R_5(7)$ and $R_5(9)$ consistent with Table 2.  Obviously there are further implied improvements in $R_r(7)$ and $R_r(9)$ for higher values of $r$.

The factors $g_k$ in Table 3 that are highlighted in bold indicate lower bounds for the limiting growth rate of $R_r(k)$ as $r$ increases,  i.e. $\lim_{\substack{r \rightarrow \infty}} R{_r}(k)^{1/r} \geq g_k$ for relevant $r$, $k$.  

\begin{table}[!ht]
  \begin{center}
    \includegraphics{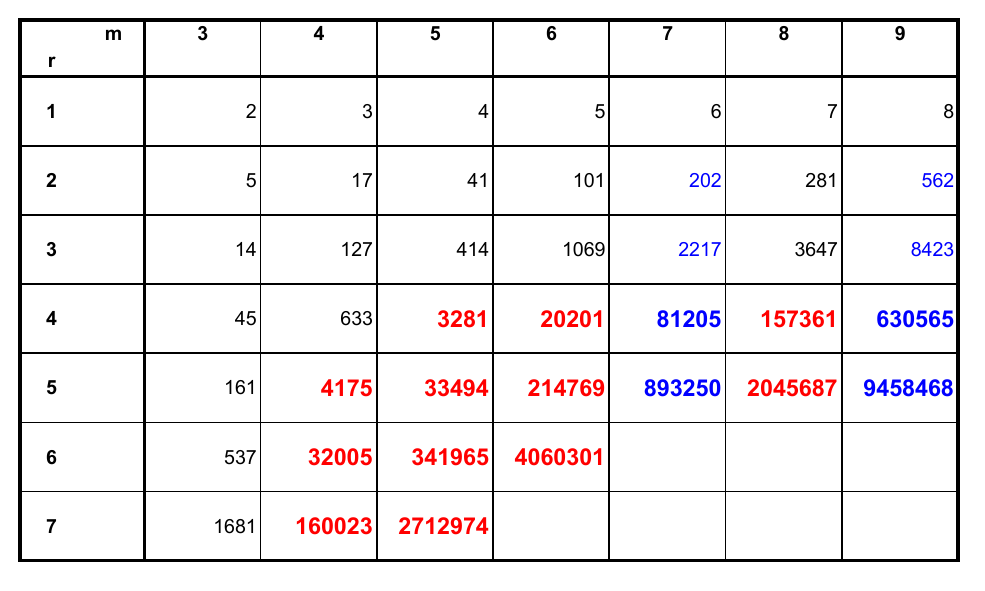}
  \end{center}
  \caption{\label{Table:02} Highest order of linear Ramsey graphs known to the author (updated).}
\end{table}

\medskip

\begin{table}[!ht]
  \begin{center}
    \includegraphics{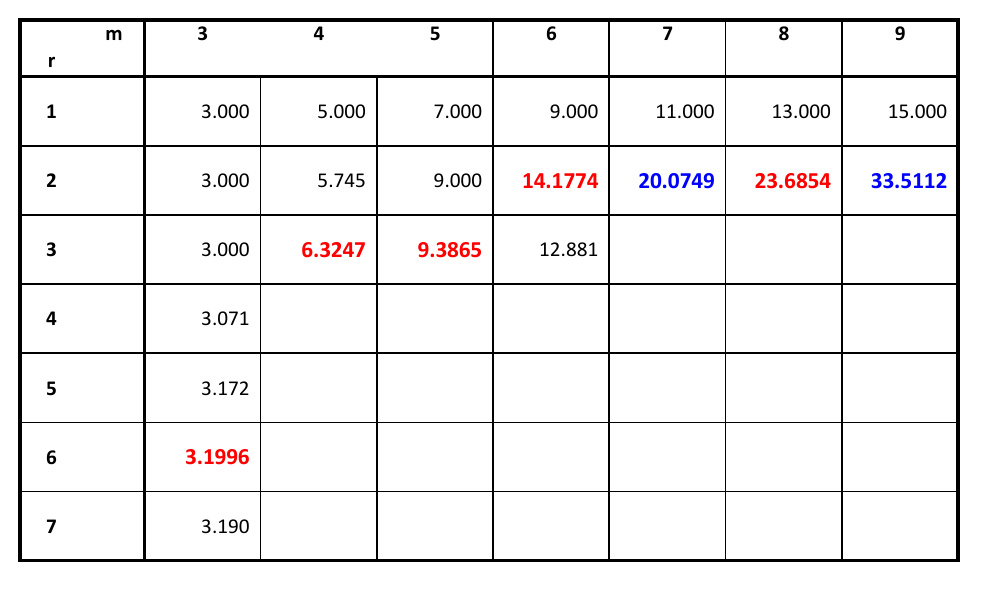}
  \end{center}
  \caption{\label{Table:03} Factors $g_k$ calculated from the data in Table 2.}
\end{table}

\FloatBarrier

\subsection*{Acknowledgements}
I once again record my warmest thanks to my wife Joan for her continued support for this work.


\begin{thebibliography}{6}

\bibitem{FrSw} H.~Fredricksen and M.M.~Sweet, \newblock Symmetric Sum-Free Partitions and
Lower Bounds for Schur Numbers, \newblock {\em Electron. J. Combin.}, {\bf 7} (2000),  {\#}R32.

\bibitem{Mathon} R.~Mathon, \newblock Lower Bounds for Ramsey Numbers and Association Schemes, \newblock {\em Journal of Combinatorial Theory}, Series B, {\bf 42}, (1987), 122--127.

\bibitem{RadzDS} S.P.~Radziszowski, \newblock Small Ramsey Numbers, \newblock {\em Electron. J. Combin.} (2017), {\#}DS1.

\bibitem{Rowley1} F.~Rowley, \newblock Constructive Lower Bounds for Ramsey Numbers from Linear Graphs, \\ \newblock {\em Australasian J. Combin.} {\bf68(3)} (2017), 385--395.  

\bibitem{XSR2} X.~Xu, Z.~Shao and S.P.~Radziszowski, \newblock  More Constructive Lower Bounds on Classical Ramsey Numbers, \newblock {\em SIAM J. Discrete Math.}, {\bf 25} No.1 (2011), 394--400.  

\bibitem{XXER} X.~Xu, Z.~Xie, G.~Exoo and S.P.~Radziszowski, \newblock Constructive Lower Bounds on Classical Multicolour Ramsey Numbers, \newblock {\em Electron. J. Combin.} {\bf 11} (2004), {\#}R35.

\end{thebibliography}

\end{document}